\documentclass{article}

\usepackage[utf8]{inputenc}
\usepackage[T1]{fontenc}
\usepackage{lmodern}

\usepackage[nopatch=footnote]{microtype}

\usepackage[
square,
authoryear,
]{natbib}

\setcitestyle{aysep={}}

\usepackage{amsmath}
\usepackage{amsfonts}
\usepackage{amssymb} 

\usepackage[amsthm]{ntheorem}

\newtheorem*{theorem-non}{Theorem}
\newtheorem*{corollary-non}{Corollary}

\usepackage[dutch,french,ngerman,greek.polytonic,english]{babel}
\usepackage{csquotes}

\usepackage{bussproofs}

\usepackage{enumitem}

\usepackage{tikz}
\usetikzlibrary{
babel,
cd,
}

\usepackage{hyperref}

\newcommand{\dutch}[1]{\foreignlanguage{dutch}{#1}}
\newcommand{\french}[1]{\foreignlanguage{french}{#1}}
\newcommand{\german}[1]{\foreignlanguage{ngerman}{#1}}

\newcommand{\weg}[1]{}

\newcommand*{\from}{\colon} 
\newcommand{\numberset}[1]{\mathbb{#1}}
\newcommand{\numnat}{\numberset{N}}

\DeclareMathOperator{\lth}{\ell h}

\newcommand{\pair}[2]{\mathopen{\langle}\,#1, #2\,\mathclose{\rangle}}

\usepackage{etoolbox}

\title{On distance and proximity between Dummett and Brouwer}

\author{Mark van Atten\thanks{Husserl Archive (CNRS\,/\,ENS), 45 rue d'Ulm, 75005 Paris, France. Email: \mbox{vanattenmark@gmail.com}}}

\date{March 11, 2026}

\begin{document}

\maketitle

\begin{abstract}
This paper asks what Brouwer might have replied to Dummett's interpretation of intuitionism.
Complementing earlier literature,
it treats
Dummett's rejection of the ontological approach;
the charge of psychologism and solipsism; 
indefinite extensibility; 
and predicativity.
It is argued that Dummett's direct arguments against Brouwerian intuitionism
do not settle the matter,
and that, 
on the latter two themes,
Dummett's position comes closer to Brouwer's than his own account suggests.
The remaining philosophical distance,
however, 
is substantial.
\end{abstract}

\setcounter{tocdepth}{1}

\tableofcontents

\section{Introduction}

While Dummett in developing his anti-realism
in logic and mathematics
took his cue from intuitionism as proposed and developed by Brouwer and Heyting,
it was clear from the outset
that there were differences.%
\footnote{John Crossley, 
who first met Dummett in 1960 at the latest, 
and I have been trying to find out more details about the history of Dummett's engagement with intuitionism. 
But apart from Prof.~Crossley's precious memories and Dummett's short remarks in his 
interview with Schulte
\citep[appendix, pp.~175–176]{Dummett1993c}
and in his
intellectual autobiography
\citep[pp.~15–16]{Dummett2007a},
so far we have found few sources.}
These concern
not only the philosophical bases,
but also
the mathematics that can be founded on them.
Much of the discussion of Dummett's 
project in the philosophy of logic and mathematics has 
focused on the question of how realists may respond.
But there is also
the question:
What might Brouwer have replied?

Various aspects and instances of that question
have led to a smaller but significant literature.
Among the first to contribute was
Göran Sundholm,
who 
in his review of 
the first edition of
\emph{Elements of Intuitionism}
wrote:
\begin{quote}
Furthermore it should be remarked that it is not clear that the 
argument based on the proof-theoretical meaning theory, with its 
presupposition of normalisation and harmony, actually justifies 
intuitionistic \emph{analysis}. 
The reviewer is willing to accept that systems 
adequate for a formalisation of Bishop's Constructive mathematics, 
e.g., Martin-Löf's type theory, 
can be justified using the Dummett 
argument. 
He cannot see, however, how to proceed 
to the intuitionistic notion of choice sequence
 within the Dummett framework. 
\citep[p.~94]{Sundholm1979}
\end{quote}
And some 15 years later,
in a paper written for Dummett,
\enquote{Vestiges of realism}:
\begin{quote}
In my opinion, 
\enquote{constructivism} 
deserves by
far the preference above 
\enquote{intuitionism} 
as an appellation for the sort of
mathematics to which Prof.~Dummett is attracted. 
Indeed, Heyting
gives no meaning-explanations 
for typically intuitionist notions such as
choice sequences or propositions involving the creative subject and, as
Troelstra observes, it is uncertain, whether the Brouwerian notion of
choice sequence permits a molecular meaning theory of the sort envisaged by Prof.~Dummett.
Furthermore, the intuitionist doctrine of
the languagelessness of mathematics seems very hard to square philosophically with Prof.~Dummett's meaning-theoretical concerns.
\citep[p.~154]{Sundholm1994}
\end{quote}
Unfortunately,
in his \enquote{Reply to Sundholm} \citeyearpar{Dummett1994},
Dummett  does not address the points made here.

Of prior
comparative discussions on Dummett and Brouwer
I mention the following.
\citet[section 1 and p.~227]{Parsons1986}
discusses the fact that Brouwer's concept of intuition is virtually absent from Dummett's writings. 
\citet{Tsui-James1999}
argues that,
at a largely programmatic level,
Dummett and Brouwer share
a naturalism that is not scientistic but consists in recognising a fundamental theoretical role for properties of the (thinking, speaking)
 subject,
and an approach that is metaphysically general in that their foundations of mathematics 
are based on considerations that are not specific to mathematics, 
but concern thought about objects in general.
The main point of 
\citet{Tieszen2000} 
is that an understanding of the intuitionism of Brouwer and Heyting
in terms of intentionality
comes closer to their views than Dummett's do,
based as the latter are on the linguistic turn.
Sundholm and I have remarked on differences between
Brouwer's \enquote{canonical demonstrations}
and
Dummett's
\enquote{canonical proofs}
\citetext{\citealt{Sundholm.Atten2008};
\citealt[pp.~29–33]{Atten-E}},
a contrast arising within the broader opposition 
between ontological and meaning-theoretical approaches to intuitionism.
In 
\citet{Atten2022},
I examined Dummett’s objection to the ontological route in detail, 
and argued that his critique presupposes a conception of mathematical procedures that Brouwer would not accept. 

The latter analysis forms part of
the immediate background to the present paper.
Rather than extending that discussion in its own terms, 
my aim here is to complement it 
and the other earlier work 
by considering what Brouwer might have replied on three further topics.
Section 2 briefly recapitulates my earlier analysis of
Dummett's objection,
since it underlies much of his interpretation of intuitionism 
and bears directly on the comparisons pursued here.
Section 3 turns to the first topic,
that of the recurring charge against Brouwer 
(in Dummett and elsewhere)
of
psychologism and solipsism. 
Drawing on
formulations
from Brouwer's publications, notes, and correspondence,
it documents that he saw things differently.
Sections 4 and 5 
emphasise how 
Dummett's views 
on 
indefinite extensibility (the second topic)
and,
perhaps more surprisingly,
on
predicativity
(the third)
come
closer to Brouwer’s than might have been expected, even if the convergence is incomplete.
Taken together, these considerations suggest that the
proximity between Dummett and Brouwer is greater, and the
remaining distance more demanding to articulate, than the genuine opposition between
meaning-theoretical and ontological routes might suggest.%
\footnote{%
A resemblance between Dummett and Brouwer that is not,
or not immediately, 
philosophical in nature,
is that between their ways of expressing themselves.
Since these are so characteristic,
not surprisingly each has been remarked upon before.
Compare:
\begin{quote}
Some readers of Dummett would say that it was ironic that he was so
preoccupied with style, since his own prose left much to be desired.
It is true that his sentences often displayed a rather unwieldy
complexity. 
But they also displayed an acute sensitivity to the
structure of the thoughts that they were intended to convey; 
and that fact, 
combined with the precision with which Dummett chose his words,
meant that there was a real clarity about his writing, 
however lacking it might have been in facility. 
The writing was in some respects like
the man~– 
marked by honesty and integrity, though it could at times be
difficult.
\citep{Moore2011}
\end{quote}
and
\begin{quote}
Abstraction coupled with a lively geometric intuition
are the hallmarks of Brouwer's work. 
Also bearing witness 
to the impeccable exactness
of his reasoning is his style, 
the manner of formulating sharply,
which terrified entire generations,
but today seems so natural that we now find it easier to read Brouwer
than his contemporaries. 
He taught his courses and gave his lectures in the same style~– 
impressive, 
but also breathtaking.
\citep[p.~340, trl.~mine]{Freudenthal.Heyting1968}
\end{quote}
On Brouwer's style,
see also
\citet[p.~263]{Dalen2013}.}

\section{Dummett's argument against the ontological route begs the question}%
\label{L008}
The most fundamental point of divergence 
concerns Dummett’s rejection of what he calls the ontological route to intuitionistic logic. 
In his influential essay 
\enquote{The philosophical basis of intuitionistic
logic}
\citeyearpar{Dummett1978a}, 
Dummett explores two distinct paths for embracing
intuitionistic logic rather than classical logic 
within number theory.
One is meaning-theoretical, 
advocated by himself, 
and the other
ontological, 
notably associated with Brouwer and Heyting.
Dummett
maintains that the meaning-theoretical route remains viable but firmly
rejects the ontological route.

Note that the potential force of Dummett's argument,
which he appealed to ever since, 
is in no way diminished
by his ascription to Brouwer of
\enquote{the development of a theory of meaning for mathematical statements}
\citep[p.~471]{Dummett1993a}. 
The question is rather whether intuitionistic logic is argued for on the basis of such a theory of meaning,
or on the basis of a languageless ontology,
which we describe in language but in no way depends on it.%
\footnote{On
Brouwer's meaning specifications for logic,
see also
\citet[p.~32]{Atten.Sundholm2017}.}

If Dummett's rejection stands, 
the conversation would seem to end already at the level of basic methodology.
If, however, the objection begs the question, then the dialectical situation is reopened.
In
\citet{Atten2022}, 
I have argued that it does beg the question.
Let me summarise that argument.

According to the ontological perspective,
mathematical objects are mental constructions that come into existence
through certain activities of the human mind. 
A mathematical truth then
fundamentally depends 
upon the procedures through which we demonstrate
them:
mathematical truths come into being with our acts.
As Heyting put it:
\begin{quote}
If mathematics consists of mental constructions, 
then every mathematical
theorem is the expression of a result 
of a successful construction. 
The proof of
the theorem consists in this construction itself, 
and the steps of the proof are
the same as the steps of the mathematical construction. 
\citep[p.~107]{Heyting1958C}
\end{quote}

Dummett argues that this perspective is untenable.
His critique specifically targets one specific consequence of the
ontological intuitionist's standpoint: 
a \emph{decidable}
true mathematical statement is true in virtue of our actually
having effected a proof which justifies the statement.

As I reconstruct it,
Dummett's
argument consists of two sub-arguments,
which I call A and B.
Argument A seeks to establish that decidable statements have their truth values
already before a corresponding decision procedure is carried out.
Argument B leads from an observation about the structure of Argument A to the
rejection of ontological intuitionism.

Argument A runs as follows.
Consider a decidable number-theoretical predicate 
\(P\) 
and a
natural number 
\(n\).
\(P\) is associated with a decision procedure
\(D\), 
and by the nature of such procedures,
we know that 
when applying \(D\) to \(n\),
we will
inevitably determine that either 
\(P(n)\) holds or 
\(\neg P(n)\) 
holds. 
Crucially,
such a procedure requires no additional information and does not rely
upon chance or external contingencies.
As Dummett writes, 
\begin{quote}
No further circumstance could be relevant to the result of the procedure – this is
part of what is meant by calling it a computation;
and, since at each step the
outcome of the procedure is determined,
 \emph{how can we deny that the overall
outcome is determinate also}? 
If we 
\emph{yield to} 
this line of thought, 
then we must hold that every statement
formed by applying a decidable predicate to a specific natural number
already has a definite truth-value, true or false, although we may not know it.
\citep[p.~245, emphasis mine]{Dummett1978a}
\end{quote}

So Dummett infers
from
\begin{quote}
If, 
for a given \(n\), 
we were to carry out procedure \(D\), 
then
we should find that 
\(P(n)\) 
holds, 
or, 
if we were to carry out procedure \(D\),
then we should find that 
\(\neg P(n)\) 
holds. 
\end{quote}
to
\begin{quote}
For a given \(n\), 
the statement 
\(P(n)\) 
already has a determinate truth-value 
prior to the execution of the procedure \(D\).
\end{quote}

Dummett reasons 
 that
the certainty that procedure \(D\) yields a definite result
is indicative
of the determinateness of 
the truth value of the proposition 
\(P(n)\) 
before we
even begin the procedure.
This contradicts the ontological intuitionist's view
that truth values of mathematical propositions
become determinate only through our mathematical acts.

Dummett then observes,
and this begins his argument B, 
that nothing in argument A
depended on an ontological view on mathematics;
the ontological intuitionist is therefore wrong to claim
that a decidable true mathematical statement is true
in virtue of our actually having effected a proof which justifies it.

Now Dummett does not hide the fact that the
crucial inference step in argument A
is not strictly deductive
(see the phrases I emphasised above).
This is how he presents the argument later:
\begin{quote}
The intuitionist sanctions the assertion, 
for any natural number, however large,
that it is either prime or composite, 
since we have a method that will, at least in
principle, decide the question. 
But suppose that we do not, and in practice
cannot, apply that method: 
is there nevertheless a fact of the matter concerning
whether the number is prime or not? 
\emph{There is a strong impulse to say that there
must be}: for surely there must be a definite answer to the question what we
should get, 
were we to apply our decision method. 
\citep[pp.~296–297, emphasis mine]{Dummett1994}
\end{quote}
To deny this conclusion would,
in Dummett's view,
be to display a
\enquote{resolute} 
\enquote{hardheaded} 
scepticism concerning subjunctive conditionals 
\citep[p.~247]{Dummett1978a}.
Remarkably,
Dummett does not quite stop to elaborate on the fact 
that this is 
an inference to the best explanation. 
About that kind of inference,
around 1980 he had said to Timothy Williamson,
in a different context:
\begin{quote}
The difference between us about how to do philosophy is that you think that
inference to the best explanation is a legitimate method of argument in philosophy, 
and I don’t. 
\citep[Dummett, quoted in][p.~264]{Williamson2016}
\end{quote}

In reply to Dummett,
we distinguish between two types
of procedures: 
investigative and generative. 
Investigative procedures
reveal truths that already exist
prior to the execution of the procedure
itself. 
By contrast, 
generative procedures 
precisely through their execution
create or bring into being
that in virtue of which a certain proposition 
about their outcome
will become true.

While Dummett's abductive inference may 
indeed be highly plausible
for investigative procedures,
it has,
by definition,
no plausibility whatsoever
for generative procedures.
But it is in terms of the latter
that the ontological intuitionists' perspective
is framed.
For them,
the procedures carried out in our mind
do not verify or uncover
pre-existing truths; 
rather, they 
bring mathematical truth about. 
In this sense,
Dummett's argument against ontological intuitionism
is question-begging.
That is of course not to show that his conclusion is false,
but it does show that the ontological route cannot be considered closed on the grounds he advances.%
\footnote{A recent discussion of Dummett’s meaning-theoretical route, 
in the context of Kuhn’s account of scientific revolutions, 
is \citet{Castro2025}, 
who argues that intuitionistic mathematics cannot be understood as a merely semantic revolution.}

\section{Brouwer's intuitionism is not psychologistic, and in some sense not solipsistic}

Even if
the ontological route cannot be ruled out by Dummett's argument, 
a distinct but related divergence concerns his characterisation of Brouwer’s standpoint itself. 
Indeed,
Dummett has always been adamant that Brouwer's intuitionism 
is psychologistic and solipsistic,
and regarded it untenable on those counts.

\subsection{Dummett}

The following passages make that assessment explicit.
Without attempting to be exhaustive,
I quote several,
and at some length,
because they show how it was a constant factor in Dummett's thought over a long period,
and because they relate back to the difference between the meaning-theoretical and the ontological approach,
the latter being understood in Brouwer's mentalist version.

\begin{enumerate}
\item
From
\enquote{The philosophical significance of Gödel's theorem}, 1963:
\begin{quote}
What we have been considering does bear on another intuitionist thesis:
that a mathematical proof or construction is essentially a mental entity,
something that may be capable of being represented by an arrangement of
symbols on paper, but cannot be identified with it. 
\textelp{}
Intuitionist language on this matter is, rightly, repugnant to anyone who has
grasped the point of Frege's repudiation of 
\enquote{psychologism}, 
of the introduction of strictly psychological concepts into logic or mathematics. 
\citep[p.~200]{Dummett1978b}
\end{quote}

\item
From 
\emph{Frege. Philosophy of Language}:
\begin{quote}
Brouwer’s writings are steeped in psychologism:
the faith that intuitionism is a tenable philosophy
of mathematics involves the faith that it is possible
to purge it of its psychologistic form.
But,
at least,
thanks for Frege’s onslaught on psychologism,
we are able to formulate what is required.
\citep[p.~684]{Dummett1973}
\end{quote}

\item
From his review of Brouwer's \emph{Collected Works}:
\begin{quote}
 It is regrettable that, in his polemical
 writing, Brouwer never addressed himself to Frege, sometimes, indeed,
 imputing to the formalists views that would have been more justly
 ascribed to Frege. If he had, we might have had a reply to Frege's
 attacks on psychologism. At any rate, Brouwer's philosophy of 
mathematics is psychologistic through and through. The proofs in terms of
 which the meanings of mathematical statements are given are not
 identified with formal proofs; they are, rather, intuitively valid ones.
 This, of course, is entirely reasonable, and enabled the intuitionists to
 be serenely unperturbed by Gödel's incompleteness result; but Brouwer
 went further, and regarded proofs as mental constructions, often only
 imperfectly communicable by language. Frege accepted Kant's view of
 geometry as founded on a priori spatial intuition, but denied that
 arithmetic and analysis are founded on intuition at all. Brouwer, on the
 other hand, seized on the Kantian idea of our primordial intuition of
 time, from which all mathematics is generated; it is in response to the
 fundamental experience of temporal succession that the discrete sequence
 of natural numbers is constructed.
\citep[p.~609]{Dummett1980}
\end{quote}

\item\label{L006}
From
\emph{The Interpretation of Frege's Philosophy}:
\begin{quote}
Brouwer's philosophy of mathematics is, however, presented in the most 
thoroughly psychologistic manner conceivable.
Mathematical objects 
do not exist independently of us, but are products of human thought, 
whose being is to be conceived; 
and there is no attempt by Brouwer to 
distinguish reason or thought from thinking~–
mathematical objects are 
the products of our mental processes. 
Instead of presenting his philosophy 
of mathematics as a consequence of a correct theory of meaning for the 
language of mathematics, 
Brouwer regards language with scepticism, 
not because it is an imperfect instrument of communication, but because 
communication itself is possible at best in a very imperfect manner; the 
proofs whose existence supplies the only legitimate notion of mathematical 
truth are mental constructions, of which symbolic or verbal renderings 
provide a necessarily defective representation. 
His thought thus notoriously 
drives towards that most subjective form of idealism, solipsism. 
\citep[p.~66–67]{Dummett1981}
\end{quote}

\item
From
\emph{Frege. Philosophy of Mathematics}:
\begin{quote}
We are usually too impressed with the really creative ideas of Hilbert 
or of Brouwer to pay much attention to the patchy or unconvincing soil in 
which they are rooted. 
\textelp{}
We 
overlook the inadequacy of Brouwer’s repeated explanations of the genesis of 
the natural-number sequence, 
and ignore his solipsism and his failure 
to 
achieve a coherent account of the relation between mental constructions and 
their symbolic formulations.
\citep[p.~292]{Dummett1991}
\end{quote}

\item\label{L007}
Finally,
from Dummett's valedictory lecture:
\begin{quote}
Brouwer was, 
notoriously, 
a solipsist, 
or something very close to one; 
but that did not vitiate his development of a theory of meaning for mathematical statements,%
\footnote{[This does not contradict the statement in quotation \ref{L006} in this list,
\enquote{Instead of \dots};
see the remark at the beginning of section 
\ref{L008}.]} 
and a consequent revisionist programme for mathematical practice. 
The reason is precisely the flagrant untruth of his solipsism. 
Far from its being the case, 
as Brouwer maintained, 
that mathematical constructions are only imperfectly communicable, 
the very opposite is true: 
they are perfectly communicable. 
Individual mathematicians  
may have different aptitudes, 
angles of attack, 
ranges of knowledge, etc., 
but they do not have 
different viewpoints on mathematical reality: 
whatever construction one mathematician discovers, 
any other is in a position to carry out.
Just for this reason, it did not matter that Brouwer was conceiving his 
mathematical language solipsistically, 
as the analogue of a sense-datum language: 
by simply reversing the 
principle that mathematical language can only imperfectly convey mental constructions carried out by any 
one mathematician, 
it could without modification be interpreted as a language common to all 
mathematicians, 
and his theory of meaning understood in terms, 
not of individual mental constructions, 
but 
of constructions available to all. 
\citep[pp.~471–472]{Dummett1993a}%
\footnote{A passage very similar to this one up to \enquote{Just for this reason\dots} occurs later in \citet[p.~15]{Dummett2007a}.}
\end{quote}

\end{enumerate}

\subsection{Brouwer}

Brouwer has in fact made various remarks
that have a direct bearing on the characterisations Dummett gives here.
They show that at least Brouwer's intentions were different.
While it is of course a further question
to what extent those intentions can be fulfilled,
an awareness of them alters what the starting point of the discussion between Brouwer and Dummett 
on psychologism and solipsism should be.

As for solipsism,
there is indeed a fairly well known passage 
where Brouwer
argues that 
\enquote{there is no plurality of mind}:
\begin{quotation}
[C]ivilized languages, mostly being cooperative languages, 
suggest a sameness for such totally different phenomena as acts of the 
subject and acts of object individuals are. 

And this suggestion is intensified by the misleading terms civilized 
languages use to characterize the behaviour of individuals in general. 
It is 
not unreasonable to derive this behaviour from 
\enquote{reason}. 
But unreasonable 
to derive it from 
\enquote{mind}. 
For by the choice of this term the subject in its 
scientific thinking is induced to place in each individual a mind with free-will 
dependent on this individual, 
thus elevating itself to a mind of second 
order experiencing incognizable alien consciousnesses as sensations. 
Quod non est. 
And which moreover would have the consequence that the mind 
of second order would causally think about the pluralified mind of first 
order, 
then cooperatively study the science of the pluralified mind, 
and in 
consequence of this study assign a mind of second order with sensation 
of alien consciousnesses to other individuals, 
thus once more elevating 
itself, 
this time to a mind of third order. 
And so on. Usque ad infinitum. 
\citep[1239–1240]{Brouwer1949C}
\end{quotation}
Its force,
such as has,
is first of all epistemic.
Brouwer argues that 
while we can become aware of \emph{objects}
by seeing invariants in the phenomena,
we cannot become aware of \emph{other minds}
that way.
His strategy here is that of a regress argument,
leading to the absurd ascription of an experience to us that we do not have.
At the same time,
he had made it clear before that he does not think there is something
contradictory about the existence of other minds:
\begin{quote}
A very essential hypothesis in the mathematical viewing of one’s fellow-humans
is, e.g., the supposition that there is in each of them a mathematical-scientific mechanism of viewing, 
acting and reflecting,
identical to one’s own.
\citetext{\citealt[p.~48]{Brouwer1933A2}, trl.~\citealt[p.~420]{Stigt1990}, modified}%
\footnote{Also \citet[p.~598]{Brouwer1975}.}
\end{quote}

\begin{quote}
What we have said so far is all rational reflection, 
i.e., mathematical
viewing, in which the content neither of purposes nor of the objects in the world of perception
plays any part. 
It is an essential hypothesis of the understanding between
people that this rational reflection has a structure which is the same for
all individuals. 
\citetext{\citealt[p.~52]{Brouwer1933A2}, trl.~\citealt[p.~423]{Stigt1990}, modified}%
\footnote{Also \citet[p.~599]{Brouwer1975}.} 
\end{quote}

Brouwer also claims:
\begin{quote}
In default of a plurality of mind, there is no exchange of thought either. 
Thoughts are inseparably bound up with the subject. 
\textelp{}
By so-called 
exchange of thought with another being the subject only touches the outer 
wall of an automaton. This can hardly be called mutual understanding. 
\citep[p.~1240]{Brouwer1949C}

\end{quote}
Yet in whatever sense Brouwer may have been a solipsist,
it clearly did not
preclude giving sense to intersubjectivity in mathematics.
In \enquote{proposition} VIII for his thesis defence he wrote:
\begin{quote}
Human understanding is based upon the construction of common mathematical systems.
\citetext{\citealt[loose leaf]{Brouwer1907}, trl.~\citealt[p.~99]{Brouwer1975}, modified}%
\footnote{\enquote{\dutch{De verstandhouding der mensen berust op het bouwen van
gemeenschappelijke wiskundige systemen.}}}
\end{quote}
While that was written some three decades earlier,
it is a view that Brouwer held on to,
and that is behind
this statement from around 1950,
published in 1981 in an appendix to Brouwer's
\emph{Cambridge Lectures}:
\begin{quote}
The stock of mathematical entities is a real thing, for each person, and for humanity.
\citep[p.~90]{Brouwer1981}
\end{quote}

When in the dissertation itself Brouwer wrote that the language of mathematics
\begin{quote}
itself is not mathematics, 
but no more than a defective 
expedient for humans to communicate mathematics to each other and to aid their 
memory for mathematics
\citetext{\citealt[p.~169]{Brouwer1907}, trl.~\citealt[p.~92]{Brouwer1975}, modified}
\end{quote}
he expressed his principled mistrust of language in a way that at the same time acknowledges its role in mathematical practice,
both individual and collective.
In the 1920s
he was involved 
in the foundation and activity of the
\enquote{Signific Circle},
a movement aiming at the general improvement of mutual understanding
by analysing speech acts and correspondingly reforming language.%
\footnote{The texts on significs in volume 1 of Brouwer's \emph{Collected Works} are not remarked upon in Dummett's review \citeyearpar{Dummett1980}.}
Its 
\enquote{basic program},
a short text also signed by Brouwer and published as a booklet only in 1939,
states at the beginning:
\begin{quote}
The meaning of a linguistic act for the speaker and for the hearer can but partly 
be determined from the words or symbols which are employed in it, 
and it can 
only approximately be expressed in other words. 
However there is a great difference in the extent to which such acts can be dissected into and approximated by 
words. 
In the language of science 
(in particular the mathematical language) 
and (albeit to a lesser extent) 
in that of technology a fairly great stability in the meaning of words 
and linguistic acts can be attained by means of example and stepwise definition
\citetext{\citealt[p.~5]{Brouwer.Eeden.Ginneken.Mannoury1939}, trl.~\citealt[p.~447]{Brouwer1975}, modified}
\end{quote}
So while it is true that Brouwer often emphasised that language is an
imperfect means for transmitting thought, 
that is not at odds
with holding that logic and mathematics are among the areas where it fulfils its
function best.

Freudenthal, 
the intuitionistic topologist who had been Brouwer’s assistant in the 1930s, 
later remarked that 
\begin{quote}
Against the formalists Brouwer was the first to stress the communicational character of language.
\citep[p.~32]{Freudenthal1960}
\end{quote}
He wrote this in his book
\emph{Lincos},
his design for a language  for communication with extraterrestrians,
to be taught to them by showing correct and incorrect uses,
and beginning with the very basics of mathematics and logic.%
\footnote{Steen opens his review \citeyearpar{Steen1962} by declaring that
\enquote{The only reason why I do not give my logic lectures after the fashion of
Lincos is that it would take too long}.}
Brouwer was, with Heyting and Beth,
editor of the series in which it appeared,
the 
\emph{Studies in logic and the foundations of mathematics}.
After initial reservations about the book's suitability for that series,%
\footnote{Letters Brouwer to Beth of June 2
and November 8, 1959
\citep[pp.~2918, 2941]{Dalen2011b}.}
as it was no treatise on mathematical logic like the other volumes,
he wrote to Beth the unexpected but telling remark:
\begin{quote}
The manner in which the book \emph{Lincos}
has been announced by N.H.U.M.~[i.e., publisher North-Holland]
as a plea for logistics \emph{\french{à outrance}}
for increasing the chances of general understanding
has removed all my initial objections to a co-patronage
of this publication.
\citep[p.~2961, trl.~mine, French expression in the original]{Dalen2011b} 
\end{quote}
Brouwer’s name accordingly remained on the cover when the book appeared.

Overall, 
we see that
Brouwer’s critique of language concerns the limited possibility of communication itself,
as Dummett correctly diagnoses in quotation
\ref{L007}
in the list above.
Because of the uncertainty that this entails,
Brouwer considers language and communication unsuitable for the foundation of mathematics,
which is characterised by certainty.
At the same time,
Brouwer clearly sees
their role in enabling any intersubjective mathematical practice at all.

As for psychologism,
that is a position that he strongly rejected:
\begin{quote}
[t]he intention of 
[my dissertation] 
is 
not to propose considerations about which different persons can have different 
opinions, 
but to establish truths which, 
just like mathematical truths, 
anybody 
who has once understood will forever affirm. 
\citetext{\citealt[p.~328]{Brouwer1908B}, trl.~\citealt[p.~106]{Brouwer1975}}
\end{quote}

\begin{quote}
[L]ike all science of experience,
that is  to say generalisation from observation,
psychology presupposes mathematics
at least up to the first infinite cardinal number inclusive,
and thus as a matter of principle rests on the intuition
of mathematical induction.
\citep[p.~200, trl.~mine]{Brouwer1911A}%
\footnote{The translation in \citealt[p.~121]{Brouwer1975},
translates the original sentence up to,
but not including,
the phrase about the intuition of induction.
\enquote{\dutch{[Z]ooals alle ervaringswetenschap,
dat is waarnemingsgeneraliseering,
onderstelt de psychologie de wiskunde 
minstens tot en met de eerste oneindige machtigheid,
steunt dus reeds principieel op de intuitie der volledige inductie.}}}
\end{quote}

These are sources that
for a long time 
Dummett could have been expected to know,
but the first volume of the \emph{Collected Works},
which contains all just quoted,
was published in 1975,
and Dummett published his review of them,
from which I quoted above, 
in 1980.%
\footnote{On p.~609 of that review, 
he regrets that Brouwer's text from 1911,
from which the quotation above was taken,
is included only in part.}
Be that as it may,
the point here is about Brouwer's positioning,
not about Dummett's knowledge of it.
The same applies in spades to archive material
that was published much later.
In 1949,
Brouwer writes in a letter to Van Dantzig,
who had criticised a so-called Creating Subject argument
on psychologistic grounds:
\begin{quote}
My conviction that
psychologistic interpretations of intuitionistic mathematics,
however interesting,
can never be adequate,
has,
if possible,
even been strengthened by your comments.
\citetext{\citealt[p.~2460]{Dalen2011b}, trl.~\citealt[p.~439]{Dalen2011a}, modified}
\end{quote}

Incidentally,
it is not known whether Brouwer ever knew Frege's attack on psychologism.
In
1973,
Dummett in his first book on Frege 
wrote
that
\enquote{Brouwer appears
to have been totally unaware of Frege's existence}
\citep[p.~xxv]{Dummett1973}.
Given that only Brouwer's published work was available then,
that was a reasonable conclusion to draw.
But
in the meantime Brouwer's notebooks towards his 1907 dissertation
have been found,
and there Frege appears.
Brouwer remarks that Frege's criticism of Hilbert's view on definitions is correct,
with reference to 
\enquote{\german{Über die Grundlagen der Geometrie II}}
\citep{Frege1903}. 
Furthermore,
it has turned out 
that 
by 1907, 
Frege’s 
\emph{Grundlagen} 
and 
\emph{Grundgesetze} 
both
\emph{were} 
present in the library of 
the University of Amsterdam
\citep[p.~223n11]{Kuiper2004}.
(They had presumably been bought on recommendation by 
the lecturer Gerrit Mannoury,
who refers to both in his 1909 book
\emph{\german{Methodologisches und Philosophisches zur Elementar-Mathematik}}.
It goes back to his lecture course on the philosophy
of mathematics
that he taught from 1903 on and which
Brouwer attended as a student.)

\subsection{Comments}

We have seen that
Brouwer's notion of the subject who carries out 
the mental constructions was
not psychologistic,
or at least not meant to be.
Brouwer in his extant writing 
does not quite go on and elaborate on his notion of the subject.
It has been proposed that the notions
of transcendental subjectivity
as in
Kant or in Husserl
can play the required role.
Above we saw that Dummett identifies Brouwer as a Kantian,
as indeed Brouwer himself at first did.
But 
to my mind,
Brouwer came to move away from Kant and closer to a view like Husserl's.%
\footnote{For documentation, discussion, and further references,
see
\citet[ch.~6]{Atten2004},
\citet{Atten2007},
and \citet{Atten2010}.}

Now Dummett,
while deploring the psychologism that he believed
to be present in intuitionism,
has also written,
in
\emph{The Interpretation of Frege’s Philosophy},
that
\enquote{a non-psychologistic form of idealism is possible, and indeed the only version worth considering} 
\citep[p.~68]{Dummett1981}.
But for him,
moving from Kant to Husserl would not do here.
In the same book he writes that
\enquote{in his later work, Husserl 
slipped back into something rather hard to distinguish from psychologism}
\citep[p.~56]{Dummett1981}.%
\footnote{Also p.~59: \enquote{if Frege followed Husserl's 
later philosophical career, in particular if he read the \emph{Ideen} of 1913, it can 
have been no surprise to him that Husserl relapsed into something indistinguishable from transcendental idealism, 
of which he claimed phenomenonology to be the first scientific version.}}

On the other hand,
Dummett had by then also recognised that 
his own meaning-theoretical approach
requires us to accept certain
premises concerning our capacities:

\begin{quote}
What differentiates such a [verificationistic] theory from one in which truth is the central notion is, 
first, 
that meaning is not directly given in terms of the condition for a 
sentence to be true, 
but for it to be verified; 
and, secondly, that the 
notion of truth, 
when it is introduced, 
must be explained, 
in some 
manner, 
in terms of our capacity to recognize statements as true, 
and 
not in terms of a condition which transcends human capacities.
\citep[p.~116]{Dummett1976}
\end{quote}

It is clear that Dummett would not have intended
to understand such a premise about our capacities
as a psychological one.
To be sure,
the idea that Dummett would accept the existence of
such a
non-psychological premise about our capacities,
which must mean our \emph{mental} capacities,
would seem to contradict a passage like the following
 (commenting not on Brouwer but on Dedekind):
\begin{quote}
[T]here is no such thing as 
\emph{the} 
human mind, only individual minds. 
The metaphor 
[of mathematical objects as creations of the mind]
is dangerously psychologistic, 
tempting us to scrutinise 
the internal operation of our minds. 
\citep[p.~311, original emphasis]{Dummett1991}
\end{quote}
Thus there seems to be a strong tension in Dummett on this point.
For the present purpose,
I will leave it at this observation,
and turn from the question of how Brouwer’s standpoint is to be understood 
to a structural comparison in which instances of Dummett's proximity to Brouwer become explicit.

\section{Dummett's indefinitely extensible concept and Brouwer's denumerably unfinished set are similar, but not co-extensive}%
\label{L009}
The first instance of such proximity is found in the case of indefinite extensibility, 
and can be traced to a common inspiration from Russell.

\subsection{Dummett}

Dummett gave his first definition of
indefinitely extensible concept
in his 1963 paper 
\enquote{The philosophical significance of Gödel's theorem}:
\begin{quote}
A concept is indefinitely extensible if, 
for any definite characterisation of it, 
there is a natural extension of this characterisation, 
which yields a more inclusive concept; 
this extension will be made according to some general  principle 
for generating such extensions, 
and, 
typically, 
the extended characterisation will be formulated by reference 
to the previous, 
unextended, 
characterisation.
\citep[p.~195]{Dummett1978b}
\end{quote}

And this is a later one from his paper
\enquote{What is mathematics about?}:
\begin{quote}
An indefinitely extensible concept is one such that, 
if we can form a definite conception of a totality 
all of whose members fall under that concept, 
we can, 
by reference to that totality, 
characterize a larger totality all of whose members fall under it.
\citep[p.~441]{Dummett1993b}
\end{quote}

These definitions are Dummett's acknowledged
absorption of
Russell’s 
notion of
\enquote{self-reproductive processes and classes}.
Russell 
defined those in
\enquote{On some difficulties in the theory of transfinite numbers and order
types} of 1906:
\begin{quote}
\mbox{}[A]ccording to current logical assumptions,
there are what we may call self-reproductive processes and classes. 
That is, 
there are some properties such that, 
given any class of terms all having
such a property, 
we can always define a new term also having the property
in question. 
Hence we can never collect all the terms having the said
property into a whole; 
because, 
whenever we hope we have them all, 
the
collection which we have immediately proceeds to generate a new term
also having the said property.
\citep[p.36]{Russell1906}
\end{quote}

\subsection{Brouwer}

Shortly after Russell's paper appeared,
Brouwer defended his dissertation,
in February 1907.
There 
we find a very similar definition:
\begin{quote}
\mbox{}[H]ere we call a set
[\emph{\dutch{verzameling}}]
denumerably unfinished
if it has the following properties:
we can never construct in a well-defined way
more than a denumerable subset of it,
but when we have constructed such a subset,
we can immediately deduce%
\footnote{[The verb here is \enquote{afleiden}~– not in a logical sense, but obtaining one construction starting from another.]} 
from it,
following some previously defined mathematical process,
new elements which are counted to the original set.
But from a strictly mathematical point of view this set does not exist as a whole,
nor does its power exist;
however we can introduce these words here as an expression
for a known intention.%
\footnote{[\enquote{Intention} translates \enquote{\dutch{bedoeling}}.]}
\citetext{\citealt[p.148]{Brouwer1907}, trl.~\citealt[p.82]{Brouwer1975}}
\end{quote}
Brouwer's definition is motivated 
by a preceding discussion why
Cantor's second number class
is not,
as a whole, 
a constructible entity;
it is a denumerably unfinished set.
Following this definition
are
the further examples of
the whole of the definable points on the continuum
and a fortiori the whole of all possible mathematical systems.
The latter example was given this form
in Brouwer's student notebooks:
\begin{quote}
The collection of mathematical theorems
is,
among others,
also a set that is denumerable,
but never finished.
\citep[VIII-44, trl.~mine]{Brouwer19051907}%
\footnote{\enquote{\dutch{Het aantal wiskundige stellingen is o.a.~ook een Menge, die aftelbaar is, maar
nooit af.}}}
\end{quote}

While those notebooks,
and that dissertation itself,
naturally contain many references to the literature of the time,
there is no indication in them that he was aware of Russell’s publication.
What the notebooks 
\emph{do}
show is that
Brouwer made his first remark on
denumerably unfinished sets 
while reading
Russell’s
\textit{Principles of Mathematics}
of 1903.
In between various quotations and references to it,
mostly on geometry,
Brouwer writes:
\begin{quote}
One will never be able to resolve the whole mystery
of spaces and surfaces in such a way,
that there is
no hocus-pocus
to it anymore;
for the number of possible
buildings [i.e., constructions]
is denumerably unfinished,
hence not surveyable.
\citep[II–32, trl.~mine]{Brouwer19051907}%
\footnote{\enquote{\dutch{Men 
zou het heele mysterie van ruimte en vlakken nooit kunnen ophelderen
zóó, 
dat er geen hokuspokus meer aan is; 
immers het aantal mogelijke gebouwen
is aftelbaar onaf, 
dus niet te overzien.}}}
\end{quote}
Brouwer did not write dates in his notebooks,
but the first page of the first notebook refers to Hilbert's contribution to
the proceedings of the 1905 Heidelberg conference,
and the first reference to a paper of 1906 occurs on the first page of the sixth.

As the quoted remark is in the second,
it was very likely written in 1905.

\subsection{Comments}

On the relation between Dummett
and Brouwer
on this point,
Richard Kimberly Heck has made this appropriate observation:
\begin{quote}
The analogy between the notion of an indefinitely extensible concept
and the more familiar notion of potential infinity should be obvious. 
In a
sense, then, it is not surprising that Dummett should suggest that the correct logic to use when reasoning about indefinitely extensible concepts is
not classical logic, but intuitionistic logic 
(Dummett [1991], p. 319). 
What
Dummett has offered us, 
then, 
is a very different sort of argument for intuitionism, 
one that is much closer to Brouwer’s original motivations than
Dummett’s 
\enquote{semantic} argument is, 
and one that is specific to mathematics
rather than almost completely general.
\citep[p.~7]{Heck2013}
\end{quote}
We now see that the closeness to Brouwer’s original motivations
is even greater than that determined by the analogy with potential infinity as such;
it is determined by the closeness of 
indefinitely extensible concept 
to
denumerably unfinished set,
or as Brouwer would later say,
\enquote{ever-unfinished and ever-denumerable species}
\citep[p.~3]{Brouwer1981A}.

Yet a distance remains.
A first difference is this.
Dummett's definition
has been observed to be circular:
The definition of indefinitely extensible
refers to a
definite conception of a totality.
But what is that,
apart from the fact that it is clearly not an indefinitely extensible one? 
Dummett offers no independent understanding of definiteness.
This circularity has been carefully 
discussed in a paper by Crispin Wright and Stewart Shapiro 
\citeyearpar[esp.~265–266, 280]{Shapiro.Wright2006}.
Here I will limit myself to pointing out that
the situation is somewhat different for
Brouwer's
definition of denumerably unfinished sets.
The opposite of \enquote{unfinished}
is \enquote{finished}
and means that we possess a law to generate the object,
in this case a set.
While the notion of law is itself open-ended,
in order to recognise a given procedure as lawlike,
no general definition of laws would be needed to begin with.
In fact,
as pointed out by Skolem, Péter, and Heyting,
to argue that the constructive notion of law
should be defined as 
\enquote{recursive}
would be circular.
\footnote{For discussion and references, see \citet[sections 1 and 2]{Coquand2014} and \citet[p.~14]{Sundholm2014}.} 
Lawlikeness is a fundamental concept that we can pick up from examples.
And one of the points about Dummett's definiteness
that Wright and Shapiro make
is precisely that his definiteness is 
\enquote{too sophisticated to allow of explanation only by examples}
\citep[p.~266]{Shapiro.Wright2006}.

A second difference is that Dummett's concept
lends itself to an attempt to argue 
that the concept of natural number
is indefinitely extensible
in the following way:
\begin{quote}
We have a strong conviction that we do have a clear grasp of the totality of 
natural numbers; but what we actually grasp with such clarity is the principle 
of extension by which, given any natural number, we can immediately cite one 
greater than it by 1. A concept whose extension is intrinsically infinite is thus 
a particular case of an indefinitely extensible one.
\citep[p.~318]{Dummett1991}
\end{quote}
And:
\begin{quote}
This, however, is again to assume that we have a grasp of the totality of natural 
numbers: but do we? 
Certainly we have a clear grasp of the step from any natural number to its successor: 
but this is merely the essential principle of extension. 
The totality of natural numbers contains what, 
from our 
standpoint, 
are enormous numbers, 
and yet others relatively to which those are minute, 
and so on 
indefinitely; 
do we really have a grasp of such a totality? 
\citep[p.~442]{Dummett1993b}
\end{quote}

Brouwer may respond as follows.
On his ontological view,
what we can construct is 
not only a multiplicity of numbers,
but also,
founded on that,
one higher-order object,
namely the potentially infinite sequence
of which,
at each moment, 
the numbers we have generated so far
make up an initial segment.
We construct that higher-order object 
in an act of reflection on the acts we are carrying out.
In other words,
the claim is that this reflection makes evident
not only that after \(n\) we can construct \(n+1\),
which corresponds to Dummett's claim here;
but that it also makes evident that 
(a)~we can iterate this and that 
(b)~doing so allows us to construct, as an additional object, the one growing sequence
1, 2, 3, \ldots.
To explain this in detail,
and to make it evident at least relative to Brouwer's point
of departure based on intuition of the passage of time,
would take a phenomenological analysis
\citep{Atten2007,Atten2024a}.%
\footnote{An early example of a position where on the one hand it is recognised that certain finite sequences can always be extended,
but on the other hand it is in fact  denied that this enables one to construct one potentially infinite sequence,
is Kant's,
notably in his discussion with Rehberg
\citep{Atten2012}.}
Here
I just remark
that 
from Dummett's perspective,
Brouwerian intuitionists, 
by taking into consideration
time awareness and the mental acts it makes possible, 
import something non-mathematical into mathematics.
But from the Brouwerian perspective,
the meaning-theoretical account neglects certain discernable properties of the construction process,
and thereby a further type of evidence
that can be appealed to in making certain mathematical principles evident.

\section{Brouwer's intuitionism is predicative, for a reason that Dummett sees but does not generalise}

A second, 
and in some respects deeper, 
instance of Dummett’s proximity to Brouwer concerns predicativity, 
though only on certain occasions.

\subsection{Brouwer}%
\label{L005}

In this section I argue that
Brouwer’s intuitionism is inherently predicative.
That runs counter to a common view in the literature, 
according to which he relies on impredicative definitions
in his understanding of implication and his acceptance of
generalised inductive definitions. 
Those topics
and the reception,
or lack thereof, 
of Brouwer on this point deserve separate, more general treatment.%
\footnote{On implication, 
see the preprint \citet{Atten-E}, which will not be published in that form.
On inductive definitions,
the slides \citet{Atten-slides}.}
Here I focus on those structural features of his account that are relevant to the comparison with Dummett.
While in his writings Brouwer does not give a systematic elaboration of those features,
an unfolding of what he does say
shows that critical impredicativity cannot arise.

In Brouwer's mature intuitionism,
the theoretical notion of a collection
is that of the species;
this replaces both in name and in meaning the earlier
\enquote{set},
which we saw in section
\ref{L009}.%
\footnote{The collections called 
\enquote{spreads} 
that also figure in his later work are a
particular kind of species.}
\enquote{Mathematical species} 
are defined as
\begin{quote}\label{L183}%
properties supposable for mathematical entities previously acquired,
and satisfying the condition that,
if they hold for a certain mathematical entity,
they also hold for all mathematical entities which have been defined to be equal to it,
relations of equality having to be symmetric, reflexive and transitive;
mathematical entities previously acquired for which the property holds
are called elements of the species
\citep[p.142, original emphasis]{Brouwer1952B}
\end{quote}
He adds:
\begin{quote}\label{L184}%
With regard to this definition of species 
we have to remark firstly that, during the 
development of intuitionist mathematics, 
some species will have to be considered as 
being re-defined time and again in the same 
way, secondly that a species can very well 
be 
an 
element 
of another 
species, 
but 
never an element of itself. 
\citep[p.142]{Brouwer1952B}
\end{quote}

We see that the definition of a species \(S\)
is a form of separation,
separating certain objects out of the domain 
formed by the constructions obtained
prior to 
the act of defining \(S\).
In this sense there is an extensional component to a species.

The existence of an object therefore
precedes the existence of  any species that it may become an element of.
There may well be an impredicative characterisation of 
a certain element of a certain species,
but by the definition of species,
the object in question has already been acquired 
(i.e., constructed, or a method is known to construct it).
This is a more general consequence of the definition
than Brouwer's remark that a species can never be an element
of itself.

With its reference to the \enquote{previously acquired} objects,
and thereby to a present and past,
there is an indexical component to the definition of a species.

Now it is possible that
once we have defined a species
and continue our constructive activity,
we may come to acquire
further objects
that we show to have the property in question.
Then,
as Brouwer points out,
the species needs to be defined anew,
with the same property,
but relative to a now extended stock of acquired objects.

So in general we obtain a
sequence of species,
such that the elements of each have the relevant property,
and the elements of each are included in those of the next.
The sequence stops in case we can prove that we either have constructed
all objects with the property, or have indicated a method for doing so.
Otherwise, it will go on.

Brouwer chooses to construe such a sequence of species as one and the same species that 
\emph{grows}.
This is clear from the fact that in the quotation we saw
he speaks of
\enquote{\emph{redefining} [some species] time and again in the same way},
instead of defining a new species,
which would notably  have to have a different name.
In particular,
the 
\enquote{denumerably unfinished sets}
of Brouwer's dissertation
are now seen as growing species. 

An element of a growing species has necessarily been added to it at some stage:
\(x \in S\)
is understood as
\(\exists i{\in}\numnat . x \in S_i\).
Since species are formed by continued separation,
membership of a species is determined only relative to objects defined at prior stages.

In keeping with the fundamental intuitionistic idea that
\(\numnat\)
is not a completed totality given by comprehension,
but is given on the basis of the intuition of temporal succession,
the stage index \(i\) reflects this unfolding of time rather than a totality determined in advance.
Dummett 
\citeyearpar[p.~199]{Dummett1978b} 
gave a by now much discussed argument to the effect that the concept of natural number is impredicative, 
on the ground that it includes the proof principle of induction,
which quantifies over all collections of numbers. 
But underlying this is the idea that it is by 
defining 
\(\numnat\) 
as the collection for which the principle of induction is valid
that the natural numbers are introduced,
and that is a conception that Brouwer does not share.%
\footnote{For further analysis
of iteration and induction in Brouwer, see \citet{Atten2024a} and the references given there.}

Existential quantification,
\(\exists x{\in}S . P(x)\),
now is analysed as
\(\exists i {\in} \numnat\, \exists x {\in} S_i . P(x)\),
and 
universal quantification,
\(\forall x {\in} S . P(x)\),
as  
\(\forall i {\in} \numnat\, \forall x {\in} S_i . P(x)\).
Thus,
at no stage an element will be defined in terms of a totality to which it does not yet belong.

Of course, 
most stages \(i\) lie in the future.
My future activity may develop in a variety of ways
that I cannot survey now;
in particular, it may depend on free choices.
So 
\(x \in S 
\leftrightarrow
\exists i {\in} \numnat . x \in S_i\)
does 
\emph{not}
imply that the species 
\(S\)
is countable.
It just means that, 
on the hypothesis that at some 
past, present or future moment I have proved,
for some particular
\(x\),
the one side of this equivalence,
I can then also prove the other.
That does not amount to claiming that I now have
a construction method that generates all of \(S\).

Likewise,
for a
function
or a 
predicate
to be
defined on
\(S\)
means
to be defined and,
when necessary, 
redefined,
on each of the species
\(S_i\).
The correct application of 
\(f\)
presupposes that we have shown that its argument
is an element of \(S\),
hence of 
\(S_i\)
from certain values of
\(i\)
onward.
Since in general we have no knowledge about the future 
\(S_i\)
except that its elements have been proved to have the relevant property,
that is the knowledge on the basis of which we must show
that
\(f\)
is well-defined on the growing species.

In many contexts
the temporal parameter plays no role and is suppressed.
Exceptions are, e.g., the \enquote{Theory of the Creating Subject},
and the individuation of non-lawlike choice sequences.%
\footnote{On the former,
see \citet{Atten2018} for discussion and further references;
similarly,
on the latter,
\citet{Atten2007}.}

But it is only because of that suppression
that we can conceive of
\(S\)
as a growing species,
and that no critical impredicativity can be introduced.
It is when we forget about this suppressed parameter
that the shorthand forms can be misread as allowing
for impredicativity to arise.

In the more general constructive context,
the importance of potentially suppressed parameters was highlighted
by Kreisel:
\begin{quote}
Finally, isolation of primitive concepts, in terms of which the
others can be defined, and laws (axioms) for these primitives.
Current candidates are construction
(function) and the application
operation with proof as a suppressed parameter.
\citep[p.~121]{Kreisel1965}
\end{quote}
with the footnote
\begin{quote}
As ordinal and order of the cumulative type theory are suppressed in the
practice of set theory.
The occurrence of such hidden parameters seems essential 
in work that aims to give an analysis of informal mathematics.
\citep[p.~121n13]{Kreisel1965}
\end{quote}
In the Brouwerian case,
not only proof is an often suppressed parameter,
but so is time.
In his generative view on procedures
(section 2)
that comes naturally.

\subsection{Dummett}%
\label{L004}

On two occasions,
Dummett comes close to Brouwer's view on species as,
in general,
growing.

The first comes,
not surprisingly,
in a discussion of indefinitely extensible concepts,
in his book on Frege's philosophy of mathematics:
\begin{quote}
Better than describing the intuitive concept of ordinal number as having a 
hazy extension is to describe it as having an increasing sequence of extensions: 
what is hazy is the length of the sequence, which vanishes in the indiscernible 
distance. 
The intuitive concept of ordinal number, 
like those of cardinal 
number and of set, 
is an indefinitely extensible one.
\citep[pp.~316-317]{Dummett1991}
\end{quote}

While Dummett is here discussing the classical context
\citep[p.~315]{Dummett1991},
he would have every reason to make the same judgement
in the constructive case.%
\footnote{\citet{Wrigley2026} discusses intensional and extensional readings of indefinite extensibility,
and argues that neither will do for Dummett's argument against Platonism and classical logic.
In the present intuitionistic context, there can only be intensional readings.}
And in the more specifically intuitionistic case,
the proposal to describe certain concepts as having
an increasing sequence of extensions
is almost the same as describing them as having
one extension that grows.
For intuitionistic ordinals,
the \enquote{haziness}
does not find its expression in the length of the sequence,
because that will be
\(\omega\);
instead,
it finds its expression in the fact that the ordinals cannot be generated in one lawlike sequence.%
\footnote{For the case of ordinals,
note that since a potentially infinite sequence
is also a  legitimate construction for Brouwer,
and so, 
once begun,
will count among the 
\enquote{previously acquired}
objects,
 as part of an \(\omega\)-sequence of acts
we can perfectly well  construct order types greater than \(\omega\).
For example: In the first act,
construct the potentially infinite sequence
\(\omega\);
in the second,
construct the number 1;
in the third,
take the ordered sum
of the objects constructed at the two preceding stages,
and thereby obtain
\(\omega + 1\).}

The second place where Dummett comes close to Brouwer
is in a remark on definitions.
It is made in a discussion of Troelstra's Paradox
in the Theory of the Creating Subject.
Although that Theory is often treated
as a more peripheral topic in intuitionism than I think it is,
it is there that Dummett comes to suggest:
\begin{quote}
If we view mathematical
constructions as being effected in time,
then this must apply not only to proofs
but to definitions;
and,
in the case of an inductive definition,
at any given temporal stage the definition may have been effected only for a
part of the domain.
\citep[pp.~241]{Dummett2000b}%
\footnote{Identical in the first edition, \citet[p.~347]{Dummett1977}.}
\end{quote}
This is precisely the insight that animates Brouwer's notion of species.

But Dummett's definition of species in fact
differs from Brouwer's,
in both editions of \emph{Elements}:
\begin{quote}
[W]e must introduce the 
notion of a species, the intuitionistic analogue of the classical notion of a set, 
or, more exactly, that of a property. 
Given any well-defined domain for variables 
of quantification, for instance the natural numbers, 
we regard a species of elements 
of that domain as determined by a definite condition which any element 
must satisfy to be a member of that species; 
a condition is definite when we 
know what to count as a proof, for any element of the domain, 
that it satisfies 
that condition. 
\citep[p.~26]{Dummett2000b} 
\end{quote}
Note that
the aspect of redefinition that figures in Brouwer’s definition
is here absent.
But it is precisely that mechanism that implements Brouwer's idea that species are,
in general, growing objects.

Dummett observes that impredicative definitions of elements of a species are to be ruled out,
on pain of contradiction:
\begin{quote}
It is clear, however, that we can have no constructive justification for the full 
impredicative comprehension axiom schema, asserting the existence of a species 
of objects satisfying any statable condition, including one involving quantification 
over species of the same or higher type; for such a condition would not be 
definite unless the domain of the species-variables was already determinate, and, 
in attempting to specify that domain in terms of definite conditions for membership of a species, 
we should be committing the fallacy prohibited by Russell’s
vicious circle principle. 
It is equally clear that the predicative comprehension 
axiom schema is acceptable: there will be a species determined by any condition 
involving quantification only over already specified domains.
\citep[pp.~26–27]{Dummett2000b}
\end{quote}
But unlike Brouwer,
he does not offer a conception of species 
on which critical impredicativity is evidently not at issue;
a conception that is such that predicativity falls out of it.

Dummett exploits his perceptive remark on definitions only in that one place,
to show that from weaker assumptions than figure in Troelstra's Paradox,
paradox still ensues.
That paradox is presented and analysed in the appendix to the present paper.
Dummett's comment on his paradox,
quoted there,
diagnoses that it arises from the impredicative construal of an operator.
That is correct, of course.
But what the analysis makes explicit is that, 
in this configuration, 
elements of a species are thereby defined in a way that quantifies over a collection that cannot be determined independently of those same elements. 
It is difficult to see how to reconcile this with the idea of temporal unfolding of definitions that Dummett just formulated,
according to which definitions are effected only progressively and may at any stage be incomplete.
In fact,
as we just saw,
he had already ruled out 
impredicative definition of species
in his discussion 
of the full comprehension principle
much earlier in the book.
But that impredicative principle was ruled out on the ground of a consequence
(specifically, a contradiction),
not by citing a generative nature of species,
however he would choose to understand that.
For Brouwer,
by contrast,
species have precisely such a nature,
as they consist in a possible constructional development
on the basis of the intuition of time.
That is how he can classify their introduction as part of
the \enquote{second act of intuitionism}
(e.g., \citealt[p.~142]{Brouwer1952B}).

\section{Concluding remark}

I began by raising the question of what Brouwer,
upon reading Dummett on intuitionism,
 might have replied. 
The discussion suggests that
their conversation 
cannot be regarded as closed
at the methodological level. 
Perhaps Brouwer and Dummett would find at least some common ground
in a sufficiently general view on non-psychologistic idealism.
And in their different ways, 
both recognised the fundamental phenomenon of extensibility; 
but Brouwer’s treatment of it allows for a theory that avoids some of the complications that arise for Dummett. 
Finally, 
if Dummett were to extend his view 
on the temporal unfolding of inductive definitions
to definitions of collections more generally, 
the result would approximate Brouwer’s predicative theory of species.

Much as these considerations bring Brouwer and Dummett into closer alignment,
either dialectically or doctrinally,
Dummett’s conviction that philosophy of mathematics 
must proceed through an account of linguistic meaning 
would likely have prevented him from accepting the ontological orientation of Brouwer’s approach.
That tension was already evident in the passages from Sundholm with which we began.

Yet consider what Dummett had to say
in a very general context
about the intended effect 
of his anti-realist arguments:
\begin{quote} 
It was in my article 
\enquote{Truth} 
that I came closest to absolutely endorsing an anti-realist 
view. 
In general, however, I have tried to avoid doing that. I 
have tried to remain agnostic between a realist and an 
anti-realist view, 
but have urged that the usual justifications 
of realism are inadequate, 
and that there is therefore a large 
problem to be resolved.
\citep[p.~192]{Dummett1993b}
\end{quote} 
Perhaps,
mutatis mutandis,
a detailed exposition of the respects in which his position
comes closer to Brouwer's
than he might have expected
would have prompted him to develop and articulate more fully his reasons
for maintaining a distance.

\paragraph{Acknowledgement.}
This is the revised and expanded text of an invited talk at the Michael Dummett Centenary Conference
in Oxford, June 29–July 1, 2025.
I am grateful to the organisers for the invitation, 
and to the audience for the subsequent discussion.
For instruction and exchanges on Dummett's philosophy over the years,
I am indebted
to Dirk van Dalen and Göran Sundholm.
For recent conversation on matters coming up in the present paper
I thank
Ryota Akiyoshi,
John Crossley,
and
Wesley Wrigley.

\section*{Appendix: Dummett's Paradox}%
\label{L002}
\addcontentsline{toc}{section}{\protect\numberline{}{Appendix: Dummett's Paradox}}

This section presents 
what in section 
\ref{L004}
I referred to as Dummett's Paradox,
and reconstructs it
in terms of a weak form of Lawvere's fixed point theorem.
While the main arguments in this paper are independent of this reconstruction,
it is offered
\begin{itemize}
\item[(a)] for its explanatory value,
in that it  makes explicit the diagonal fixed point structure that Dummett’s presentation instantiates but leaves implicit;
and 
\item[(b)] to facilitate comparison
with the many other 
paradoxes
that have been 
reconstructed using that theorem.%
\footnote{There are many in \citet{Yanofsky2003}, and a few in the specifically intuitionistic context in \cite{Atten2024b}.
This appendix can be seen as an addition to the latter.
Also the paradox by Martino and Usberti
\citeyearpar[p.~150]{Martino.Usberti2018}
is readily converted to this form.}
\end{itemize}

Dummett devised his argument 
as a
variation on Troelstra's Paradox,%
\footnote{Dummett presents the latter in
\citeyearpar[pp.~345–347]{Dummett1977} and \citeyearpar[pp.~240–241]{Dummett2000b}.
It first appeared in \citet[pp.~105–107]{Troelstra1969}.
See for further discussion 
\citet[pp.~1608–1609]{Atten2018}
and
\citet[pp.~464–468]{Atten2024b}.} 
to show that 
an impredicative construal of 
\(\vdash_n\)
leads to a paradox even if it is not assumed,
as in Troelstra's,
that at each stage the Creating Subject
draws exactly one conclusion.

Here is
Dummett's own presentation of the paradox.
(Familiarity with Dummett's 
expositions of the notion of spread
and of the operator 
\(\vdash_n\)
is assumed.%
\footnote{\citet[pp.~65–67 and 338]{Dummett1977}, or \citet[pp.~47–48 and 235–237]{Dummett2000b}.}%
)
\begin{quote}
If we view mathematical
constructions as being effected in time,
then this must apply not only to proofs
but to definitions;
and,
in the case of an inductive definition,
at any given temporal stage the definition may have been effected only for a
part of the domain.
Suppose,
now,
that we are defining inductively a spread-law \(c\) which
determines
a subspread of the full binary spread \(b\).
At any given temporal stage \(m\),
\(c\) may
have been defined only over certain finite sequences,
say those whose length
falls below some bound \(\ell + 1\);
so we may consider the relativized spread-law \(c^*_m\)
which admits a finite sequence 
\(\vec{u}\) 
just in case (i)
\(\vec{u}\) is admissible under \(b\) and
(ii)
every initial segment of \(\vec{u}\) for which \(c\) is defined at stage \(m\)
is admissible
under \(c\).
If we construe 
\enquote{\(\vdash_n\)} 
impredicatively,
we must regard as well-defined that
spread-law \(c\) such that,
where \(\vec{u}\) is admissible under \(c\),
then so is \(\vec{u}^\smallfrown 1\),
and,
further,
\(\vec{u}^\smallfrown 0\) is admissible iff,
where \(m = \ell h(\vec{u})\),
\(\vdash_m \neg \exists \alpha_{\alpha{\in}c^*_m}\exists n\,\alpha(n)=0\).
Now suppose that,
for some \(\vec{u}\) of length \(m\),
\(\vec{u}^\smallfrown 0\) were admissible under \(c\).
Then
\(\vdash_m \neg \exists \alpha_{\alpha{\in}c^*_m}\exists n\,\alpha(n)=0\),
and hence,
a fortiori,
\(\neg \exists \alpha_{\alpha{\in}c}\exists n\,\alpha(n)=0\),
which
contradicts the assumption that \(\vec{u}^\smallfrown 0\)
is admissible.
We have thus shown that a finite sequence is admissible
under \(c\) iff it consists entirely of 1s.
If we are presently
at stage \(k\),
\(c\) has been defined for all finite sequences at stage \(k\),
and moreover
we have proved that 
\(\neg \exists \alpha_{\alpha{\in}c}\exists n\,\alpha(n)=0\),
and so \(\vdash_k \neg \exists \alpha_{\alpha{\in}c}\exists n\,\alpha(n)=0\).
By
the definition of \(c\),
it follows that the finite sequence consisting of \(k\) 1s
followed
by 0 is admissible,
and we have arrived at a contradiction.
This contradiction
springs solely from the impredicative character we are
attributing to the notion
expressed by 
\enquote{\(\vdash_n\)},
for which we have not here assumed the strict
interpretation;
of course, if we construe 
\enquote{\(\vdash_n\)} 
predicatively, then 
\(c\) 
is not
properly defined at all.
\citep[pp.~241–242]{Dummett2000b}%
\footnote{Identical in the first edition, \citet[pp.~347–348]{Dummett1977}.}
\end{quote}

This passage will now be read in light of the following
\begin{theorem-non}[Lawvere’s fixed point theorem, diagonal form]\hspace{-0.5em}%
\footnote{\label{L001}%
This form is decidedly weaker than the full theorem 
\citep{Lawvere1969}
in that 
(a) it is formulated only for sets;
(b) it requires
representability 
not for arbitrary mappings from \(A\) to \(B\),
but 
for one specific $f$;
and 
(c) representability is not a global requirement on the domain of $f$ but only at the diagonal.
In particular,
the full theorem includes the hypothesis that, 
for every \(f \from A \to B\), some element of \(A\) represents it at all arguments.
That impredicativity is absent here,
and instead the relation of \(a\) to itself is included directly in the condition on \(g\).
While weakening the hypotheses strengthens the implication in a logical sense, 
it thereby weakens the theorem in a structural sense,
because
the conclusion is no longer derived from a single uniform 
representability principle, 
but from the diagonal condition in a specific instance.}
Let \(A\) and \(B\) be sets, and suppose that \(g\) is defined on all diagonal pairs \((a,a)\in A\times A\). Let
\[
\Delta : A \to A\times A, 
\qquad 
\Delta(a)=\langle a,a\rangle.
\]
Let $h : B \to B\), and define
\[
f = h \circ g \circ \Delta.
\]
Assume that 
\(g\)
is such that
diagonal representability holds for \(f\),
i.e., that for each \(a \in A\),
\[
g(a,a)=f(a).
\]
Then, for each \(a \in A\),
\[
g(a,a)=h(g(a,a)).
\]
Hence each \(g(a,a)\) is a fixed point of \(h\).
\end{theorem-non}

\begin{proof}
In the assumption, replace \(f\) by its definition.
\end{proof}

\begin{corollary-non}
If \(h\) has no fixed point in \(B\), we derive a contradiction.
\end{corollary-non}

\medskip

The reconstruction of Dummett's Paradox in these terms
proceeds in three steps.

\bigskip

1. We first define the sets and functions involved.
Set
\[
A=\numnat  
\]
and
\[
B=\{\{1\},\{0,1\}\}.
\]
The elements of \(B\) are the potential sets of choices for extending a given admissible sequence.

Define
\(g\)
by
\[
g \from \{\pair{m}{m} \mid m \in A\} \to B
\]
and
\[
g(m,m)
=
\{b \in\{0,1\}\mid
1^m{^\smallfrown}b
\text{ is admissible under }c^*_m\}.
\]
The value \(g(m,m)\) becomes determined once the
definition of the admissible sequences in \(c\) has been effected up to stage \(m\).

Define
\(h:B\to B\)
by
\[
h(\{1\})=\{0,1\},
\qquad
h(\{0,1\})=\{1\}.
\]
and
\[
f = h\circ g\circ \Delta : A\to B.
\]
Thus for each $m \in A$,
\[
f(m)=h(g(m,m)),
\]
so that
\(f\) 
inverts whether \(1^m\) may be extended by 0, as recorded by \(g(m,m)\).

\bigskip

2. Now we verify
that the function \(g\) as  defined here satisfies the representation condition  of the theorem.
Dummett’s spread law
\(c\) 
states: 
\begin{enumerate}
\item[(a)] If \(\vec{u}\) is admissible, then so is \(\vec{u}{^\smallfrown}1\).
\item[(b)] \(\vec{u}{^\smallfrown}0\) is admissible iff
\(
\vdash_m \neg \exists \alpha{\in} c^*_m \exists n\,(\alpha(n)=0)
\),
where \(m=\lth(\vec{u})\).
\end{enumerate}
Note the impredicativity in the second clause:
the operator
\(\vdash_m\)
is applied to a proposition that quantifies over sequences \(\alpha \in c^*_m\),
but membership in 
\(c^*_m\)
is itself determined by admissibility conditions formulated using 
\(\vdash_m\).

For each \(m\), the statement
\[
\vdash_m \neg \exists \alpha{\in} c^*_m \exists n\,(\alpha(n)=0)
\]
is equivalent to the statement that, under \(c^*_m\), no extension \(1^m{^\smallfrown}0\)
is admissible:
\begin{align*}
\vdash_m \neg \exists \alpha{\in} c^*_m \exists n\,(\alpha(n)=0) \quad &\text{iff} \quad g(m,m)=\{1\}.
\tag{Diag.}\\
\intertext{The defining clause for \(1^m{^\smallfrown}0\) can now be rewritten as}
1^m{^\smallfrown}0\ \text{admissible under }c^*_m
\quad &\text{iff} \quad
h(g(m,m))=\{0,1\},
\intertext{and therefore also as}
1^m{^\smallfrown}0\ \text{admissible under }c^*_m
\quad &\text{iff} \quad g(m,m)=f(m),
\tag{Diag.~repr.}
\end{align*}
because,
from left to right: 
if 0 is admissible,
\(g(m,m) = \{0,1\}\),
and then the equality right-hand side follows from the preceding rewriting of the clause
and the definition of \(f\);
and
from right to left:
if \(g(m,m)=f(m)\),
then by definition of \(f\) and the preceding rewriting of the clause,
\(g(m,m)=\{0,1\}\),
and the left-hand side follows from the definition of \(g\).
Thus the index \(m\) 
of the relativised spread \(c^*_m\) functions as the  diagonal representing element
for  \(f\)  when applied to 
\(m\).  

\bigskip

3. Finally,
we verify that Dummett's reasoning toward his conclusion
can be seen as an application of the theorem.
Suppose that we are at some arbitrary stage \(k\).
Dummett proves at the outset (stage 0) that
\(\neg \exists \alpha{\in}c \exists n\,(\alpha(n)=0)\),
hence
\[
\vdash_k \neg \exists \alpha{\in}c \exists n\,(\alpha(n)=0) 
\]
Therefore,
by (Diag.),
\[
g(k,k)=\{1\}.
\tag{Global}
\]
This is derived from the global proof about the full spread \(c\).

Since the relativised spread is part of the full one, 
we also have
\[
\vdash_k \neg \exists \alpha{\in}c^*_k \exists n\,(\alpha(n)=0),
\]
which,
together with the defining clause,
entails that \(1^k{^\smallfrown}0\) is admissible under \(c^*_k\), 
whence
\[
g(k,k)=\{0,1\}.
\tag{Local}
\]
This follows immediately from information about the relativised spread \(c^*_k\) itself
(and only mediately from the global proof).

(Global) and (Local) 
are contradictory.
In terms of the theorem,
(Diag.~repr.)
gave 
\(g(k,k)=f(k)\);
together 
with the definition of
\(f\), 
this yielded a fixed point, 
contradicting the Corollary.

This contradiction can be derived for any stage \(k\);
one suffices to show that the spread \(c\) is not correctly defined.
However,
precisely because
Dummett's way of generating 
fixed points here is uniform in this way,
it is natural to argue accordingly.

\bigskip

Dummett attributes the contradiction to the impredicativity of
\(\vdash_n\).
In the above reconstruction,
we also see impredicativity of 
\(c^*_k\),
and also in the self-indexing role of \(k\): 
the parameter \(k\) 
selects the very partial spread over which the defining condition quantifies.
But these are interdependent
definitions,
and the impredicativities are inseparable from one another.
None of them can be replaced by an equivalent predicative definition.
In other words,
in predicative construals they are not well-defined,
and hence neither is \(c\),
as Dummett points out for \(\vdash_n\).

\bibliographystyle{plainnat}

\end{document}